\documentclass[11pt]{article}

\usepackage{amsmath,amsthm,amssymb}
\usepackage{nicefrac}
\usepackage{url}

\usepackage{calc, graphicx} 

\usepackage{bbm}

{\theoremstyle{definition}
  
}

\newcommand{\CA}[1][p]{C_{A}}
\newcommand{\Cb}[1][p]{C_{b}}
\newcommand{\CX}[1][p]{C_{X}}

\let\theta\vartheta

\let\phi\varphi

\let\rho\varrho

\let\epsilon\varepsilon

\title{Appendix to ``Approximating perpetuities''}
\author{Margarete Knape and  Ralph Neininger\thanks{Email: \{knape, neiningr\}@math.uni-frankfurt.de, DOI 10.1007/s11009-012-9299-2}\\
Institute for Mathematics\\ J.W.~Goethe-University\\60054 Frankfurt a.M.\\
Germany}
\begin{document}
\date{July 25, 2012}

\maketitle

\begin{abstract}
  An algorithm for perfect simulation from the unique solution of the distributional fixed point equation $Y=_d UY + U(1-U)$ is constructed, where $Y$ and $U$ are independent and $U$ is uniformly distributed on $[0,1]$. This distribution comes up as a limit distribution in the probabilistic analysis of the Quickselect algorithm.  Our simulation algorithm  is based on coupling from the past with a multigamma coupler. It has four lines of code.
\end{abstract}

\section{Introduction}

In a probabilistic analysis of the algorithm {\tt Quickselect}  Hwang and Tsai \cite{hwts02} showed that, when applied to a uniformly random permutation of length $n$ and selecting a rank of order $o(n)$, the normalized number of key exchanges  performed by  {\tt Quickselect} converges in distribution to a limit distribution $\mu$. This limit distribution is characterized as the unique probability measure  $\mu={\cal L}(Y)$ such that 
\begin{align}Y\stackrel{d}{=} UY + U(1-U), \label{rde}
\end{align} 
where $\stackrel{d}{=}$ (also $=_d$) denotes equality in distribution and $U$ is uniformly distributed over the unit interval $[0,1]$ and independent of $Y$. 

The distribution $\mu$ was  studied  in \cite{knne08}. In particular we showed that $\mu$ has a bounded, $1/2$-H\"older continuous density, $\mu$ is supported by the unit interval $[0,1]$ and we developed a method to numerically approximate the density and the corresponding distribution function. In Remark 2.9 of \cite{knne08} we noted that this is sufficient to theoretically construct an algorithm for perfect simulation from $\mu$ based on von Neumann's rejection method along  the approach taken in Devroye \cite{de01}. While the numerical approximations yield an algorithm for perfect simulation in almost surely finite time, the convergence rates of our approximations are poor and the expected running time  is infinite. We do not expect such an algorithm to terminate within our lifetimes. 

Recently, Fill and Huber \cite{fihu10} published an algorithm for perfect simulation of  a related distribution, known as the Dickman distribution and characterized as unique solution of  the distributional fixed point equation $Y=_d UY + 1$.  This algorithm is based on coupling from the past of a Markov chain with the Dickman distribution as stationary distribution. The method makes use of  a multigamma coupler and of a dominating chain to deal with the unbounded support of the Dickman distribution. In fact Fill and Huber  develop their algorithm for a more general class of distributions, the Vervaat perpetuities. Devroye and Fawzi \cite{defa10} presented a different multigamma coupler and a different dominating chain resulting in a faster coupling from the past algorithm for the Dickman distribution. Both algorithms are also fully satisfactory from a practical point of view, millions of independent samples from the Dickman distribution can be generated within seconds. 

In this note we construct a  coupling from the past algorithm for the solution $\mu$ of (\ref{rde}).  Compared to the more difficult Dickman case we benefit from the special analytic structure of the densities $\varphi_x$ of $Ux+U(1-U)$ for $x\in[0,1]$. In particular, we have  
\begin{align}\label{dens_prop}
\inf_{x\in[0,1]} \inf_{t\in[0,1/4]} \varphi_x(t) \ge 1/2,
\end{align} 
which allows for the construction of a multigamma coupler as proposed by Murdoch and Green \cite[Section 2.1]{mugr98}. This results in a fast and  simple four-line-code algorithm. 

Note that a general method described in an unpublished extension of \cite{defa10}, see  Fawzi \cite{fa07}, can also be applied to our $\mu$: In \cite[Section 4]{fa07} it is shown that when one is able to perfectly simulate from the solution of $Y=_d AY+1$ with a random $0\le A\le 1$ this can be turned into an algorithm to simulate from the solution of $Y=_d AY+B$, whenever $B\ge 0$ is bounded. Here, $(A,B)$ is independent of $Y$. Hence, this method together with the simulation algorithm for the Dickman distribution yields as well an algorithm to simulate from $\mu$. 

For general perfect simulation algorithms for another class of perpetuities see
Devroye and James \cite{deja11}. For perfect simulation algorithms from stationary distributions
of positive Harris recurrent Markov chains see Hobert and Robert \cite{horo04}.

In the field of exact simulation from nonuniform distributions it is customary  to assume that a sequence of independent and identically, uniformly on $[0,1]$ distributed random variables is available and that elementary operations of and between real number such as $+$, $-$, $/$, $*$, $\sqrt{x}$, $\log x$, etc., can be performed with absolute precision, see Devroye \cite{de86} for a comprehensive account  on nonuniform random number generation.

\section{Markov chain and multigamma coupler} 
An underlying ergodic Markov chain $(X_j)$ on $[0,1]$ having $\mu$ as stationary distribution is given as follows: For all $x\in [0,1]$, given $X_j=x$, we define   $X_{j+1}$ to be distributed as $Ux+U(1-U)$ with a uniform $[0,1]$ random variable $U$. In the context of coupling from the past a realization of such a Markov chain is usually constructed with a deterministic update function $\Phi:  [0,1]\times [0,1] \to [0,1]$ such that $X_{j+1}:=\Phi(X_j,U_{j+1})$ yields a realization of the chain, where $(U_j)$ is a sequence of independent and uniform $[0,1]$ random variables. 
A trivial choice for $\Phi$ is $(x,u)\mapsto ux+u(1-u)$. However, to make coupling of the chains   possible, we follow the construction of a multigamma coupler as described by Murdoch and Green \cite{mugr98}. 

The construction is as follows: Assume that a probability density $f$ is written as $f=f_1+f_2$ with measurable, nonnegative functions $f_1, f_2$ such that $\|f_1\|_1:=\int f_1(x) \, dx$,  $\|f_2\|_1>0$. Assume that $Y_1$, $Y_2$ are random variables with densities $f_1/\|f_1\|_1$ and $f_2/\|f_2\|_1$ respectively and that $B$ is a Bernoulli$(\|f_1\|_1)$ random variable independent of $(Y_1,Y_2)$. Then the random variable $BY_1+(1-B)Y_2$ has density $f$. 

The aim now is to obtain for the densities $\varphi_x$ of  $Ux+U(1-U)$ representations $\varphi_x= r + g_x$ as above, where $r$ is independent of $x\in[0,1]$. Typically this may not be possible since one may have $\inf \varphi_x=0$  such that a non-zero $r$ independent of $x$ does not exist. However, in our particular situation we have (\ref{dens_prop}), hence we are able to choose, e.g.,  
\begin{align}\label{def_r}
r(t):= \frac{1}{2}{\bf 1}_{[0,1/4)}(t), \qquad t\in[0,1].
\end{align}
Clearly, $U/4$ has density $r/\|r\|_1$ and let us assume for the moment that a random variable $Y_x$ with density $g_x/\|g_x\|_1$ can be simulated  via its inverse distribution function (quantile function) $G_x^{-1}$, i.e.,  ${\cal L}(Y_x)={\cal L}(G_x^{-1}(U))$. Then, with a Bernoulli$(\|r\|_1)$ random variable $B$, independent of $U$, we have that for all $x\in[0,1]$
\begin{align*}
Ux+U(1-U) \stackrel{d}{=} \frac{BU}{4} + (1-B) G_x^{-1}(U). 
\end{align*}
Hence, our update function is $\Phi^{'}:[0,1]\times \{0,1\}\times [0,1] \to [0,1]$,
 $(x,b,u)\mapsto bu/4 + (1-b) G_x^{-1}(u)$. If we construct our Markov 
chain from the past using $\Phi^{'}$, in each step there is a probability of  
$\|r\|_1=1/8$ that all chains couple simultaneously. In other words, we can just start at a 
Geometric$(1/8)$ distributed time $N$ in the past, the first instant of $\{B=1\}$ when moving back into the past. At this time $-N$ we couple 
all chains via $X_{-N}:=U_{-N}/4$ and let the chain run from there until time 
$0$ using the updates $G_{X_{j}}^{-1}(U_{j+1})$ for $j=-N,\ldots, -1$. It is shown in 
\cite[Section 2.1]{mugr98} that this is a valid implementation of the coupling from the 
past algorithm in general. 

Hence, we need to derive expressions for the functions $G_x^{-1}$ containing only
elementary operations.
It was calculated in \cite[equation (28)]{knne08} that, for all $t\in[0,1]$ we have
\begin{align*}
\varphi_x(t)= \left((1+x)^2-4t\right)^{-1/2}\left({\bf 1}_{[0,x)}(t)+2\cdot{\bf 1}_{[x,b_x)}(t)\right)
\end{align*}
with $b_x:=((1+x)/2)^2$.
Hence, with $r$ given in (\ref{def_r})
we have $\varphi_x(t)\ge r(t)$ for all $x,t \in [0,1]$. 
Note that coupling occurs faster when the function $r$ can be chosen larger. For our densities  $\varphi_x$ we could as well choose
\begin{align*}
r^\ast(t)= \frac{1}{2\sqrt{1-t}}{\bf 1}_{[0,1/4)}(t), \qquad t\in[0,1].
\end{align*}
Then we  have $\varphi_x(t)\ge r^\ast(t) \ge r(t)$ for all $x,t \in [0,1]$. However,  the subsequent  inversion of distribution functions  can be done elementary with our choice of $r$.

We need to invert the distribution  functions $G_x:[0,1]\to[0,1]$ corresponding to the normalized versions of $g_x=\varphi_x-r$. We have
\begin{align*}
G_x(y)= &\frac{1}{1-\|r\|_1} \int_0^y \varphi_x(t) -r(t) \, dt= \frac{8}{7} \left(F_x(y)-\frac{1}{2}(y\wedge \nicefrac{1}{4})\right),
\end{align*}
where 
\begin{align*}
F_x(y) := \left\{ \begin{array}{cl}
\frac{1}{2}\left(1+x-\sqrt{(1+x)^2-4y}\right), & 0\le y< x,\\
1-\sqrt{(1+x)^2-4y}, & x\le y<b_x,\\
1, & b_x \le y \le 1,
\end{array}
\right.
\end{align*}
is the distribution function of $Ux+U(1-U)$. 

The inversion of $G_x$ can be done by explicit calculations and yields 
\begin{align*}  
G_x^{-1}(z) = \left\{ \begin{array}{cl}
-\frac{7}{4}z+ \sqrt{ 7z+(1-x)^2}+x-1, &\mbox{if } x\in[0,\nicefrac{1}{4}], z \in [0, q_x],\\ [2mm]
-\frac{7}{4}z+ 2\sqrt{ 7z+9 +x(x+2)}-6, &\mbox{if } x\in[0,\nicefrac{1}{4}], z \in (q_x,r_x],\\ [2mm]
\frac{1}{256}(15+8x-7z)(1+8x+7z), &\mbox{if } x\in[0,\nicefrac{1}{4}], z \in (r_x,1],\\ [2mm]
-\frac{7}{4}z+ \sqrt{ 7z+(1-x)^2}+x-1, &\mbox{if } x\in (\nicefrac{1}{4},1], z \in [0, s_x],\\ [2mm]
\frac{1}{64} (7 + 8x - 7z) (1 + 7z), &\mbox{if }x\in (\nicefrac{1}{4},1], z \in (s_x,t_x],\\ [2mm]
\frac{1}{256}(15+8x-7z)(1+8x+7z), &\mbox{if }  x\in (\nicefrac{1}{4},1],  z \in (t_x,1],\\
\end{array}
\right.
\end{align*}
where
\begin{align*}
q_x&:= \frac{4}{7}x,\quad r_x:=1-\frac{8}{7}\sqrt{x(x+2)},\\
s_x&:=\frac{1}{7}\left(3+4x-4\sqrt{x(x+2)}\right), \quad t_x:= \frac{1}{7}(8x-1).
\end{align*}

\section{The algorithm}

Our algorithm {\tt Simulate[$Y=_dUY+U(1-U)$]} has the  form  discussed in the previous section: It draws back to a sequence of independent uniform$[0,1]$ random variables $(U_{-n})_{n\ge 0}$ and an independent geometrically distributed random variable. (This clearly can be simulated on the basis of independent uniform$[0,1]$ random variables as well.)\\

\noindent


\noindent
{\tt Simulate[$Y=_dUY+U(1-U)$]}:\\
\noindent\rule[1ex]{\textwidth}{1pt}\\

\noindent
\phantom{aaaa}$N \leftarrow \mathrm{Geometric}(1/8)$\\
\phantom{aaaa}$X \leftarrow U_{-N}/4$\\
\phantom{aaaa}{\bf for} $j$ {\bf from} $-N+1$ {\bf to} $0$ {\bf do} $X \leftarrow G_X^{-1}(U_j)$\\
\phantom{aaaa}{\bf return}$(X)$\\

\noindent\rule[1ex]{\textwidth}{1pt}
\pagebreak

\noindent 
The analysis of the complexity of this algorithm is trivial as the loop is iterated a random $\mathrm{Geometric}(1/8)$ number of times, hence, e.g., on average eight times.

In Figure 1 the histogram (normalized to area $1$) of the values of 10 million independent samples generated with {\tt Simulate[$Y=_dUY+U(1-U)$]} is plotted. This simulation was done within a few seconds. A numerical approximation of the density of $\mu$  has already been presented in \cite[Figure 1]{knne08}. 

\begin{figure}\label{bild}
	    \centering
	    \includegraphics[width=9cm]{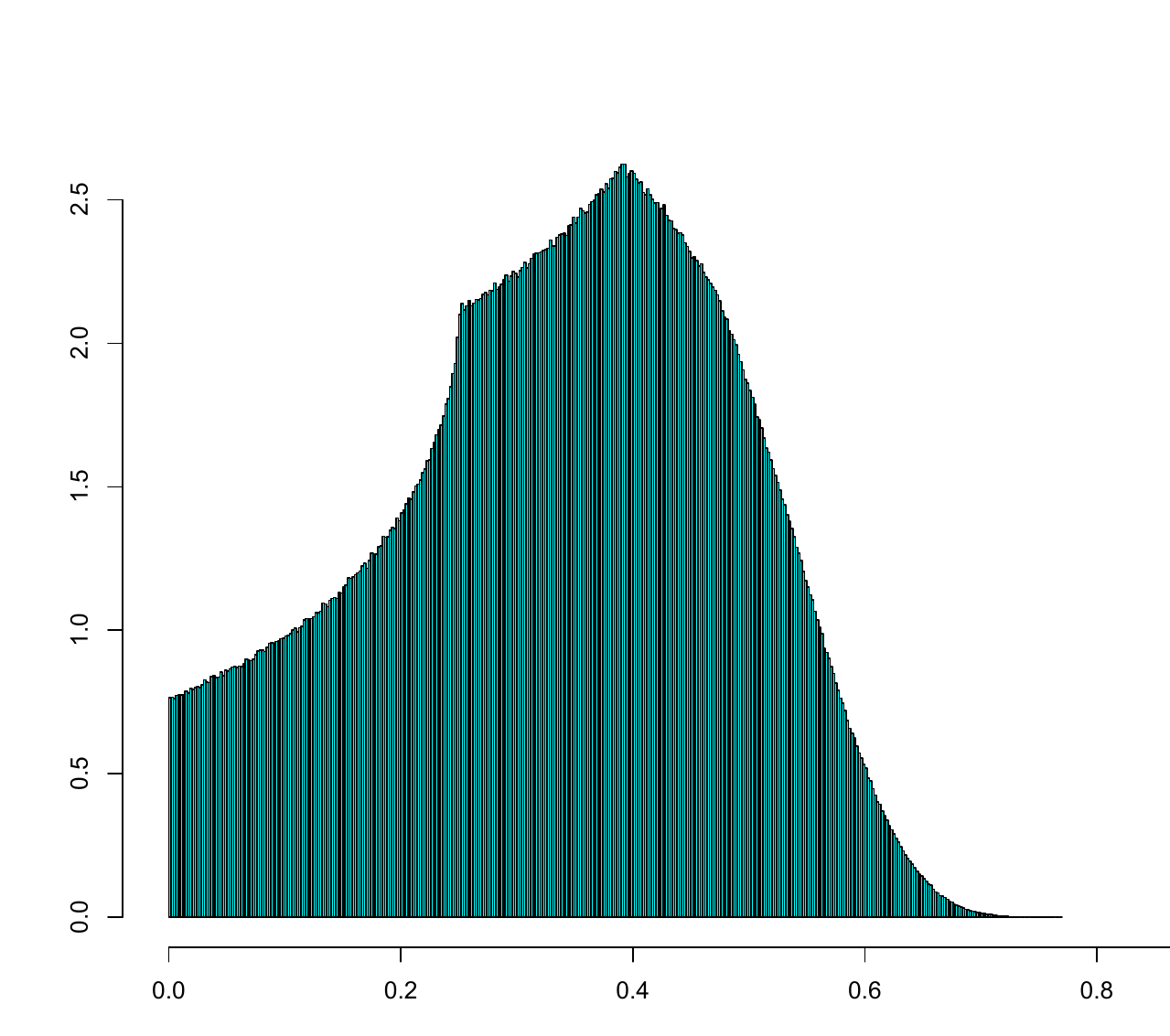}
	    \caption{Histogram of the values of 10 million independent samples from $\mu$ generated with the algorithm {\tt Simulate[$Y=_dUY+U(1-U)$]}.}
\end{figure}



\end{document}